\documentclass{article}
\usepackage{amsmath,amssymb}

\begin{document}
\newcommand{\de}{\mbox{$\delta $}}
\newcommand{\e}{\mbox{${\bf \epsilon \rm}$}}
\newcommand{\qed}{\mbox{$ \quad  \qquad \qquad \qquad \qquad \qquad \qquad \qquad \qquad \qquad \qquad
 \qquad \qquad \qquad \qquad \quad \qquad \square $ }}
\newtheorem{cor}{Corollary}[section]
\newtheorem{alg}{Algorithm}[section]
\newtheorem{lemma}{Lemma}[section]
\newtheorem{theo}{Theorem}[section]
\newtheorem{defi}{Definition}[section]
\newtheorem{exa}{Example} [section]
\newtheorem{pro}{Proposition}[section]
\newtheorem{rem}{Remark}[section]

\title{AN ALGORITHM FOR DETECTING " LINEAR " SOLUTIONS OF NONLINEAR POLYNOMIAL
DIFFERENTIAL EQUATIONS.}
\date{{\small 29 March 2006}}
\author{Stelios Kotsios\thanks{%
Faculty of Economics, Department of Mathematics and Computer Science, University of Athens, Pesmazoglou 8, Athens 10559, Greece}}
\maketitle

\begin{abstract}
A symbolic computational 
algorithm which detects " linear "` solutions of nonlinear polynomial
differential equations of single functions, is developed throughout
this paper. 
\end{abstract}

\section{Introduction}

The problem of obtaining general solutions of differential equations via symbolic algorithms
has been studied in the past by many authors, \cite{kn:kapla},\cite{kn:cheng1},\cite{kn:ritt},\cite{kn:cheng2},\cite{kn:glad1},\cite{kn:cheb}. 
These algorithms
allowed new calculation techniques to be accomplished, much more efficiently,
faster and without approximation errors.
\par \noindent
In this paper we treat with differential equations of the form:
\begin{equation}\label{arxh}
 p(x,y(x),y'(x),y''(x),\ldots,y^{(n)}(x))=0
 \end{equation}
where $p$ is a polynomial function and $y(x)$ a real 
function of a single
 variable. Our aim is to discover possible " linear " solutions
 of (\ref{arxh}). By the term "linear" we mean solutions which 
 can be obtained by solving linear differential equations. Our 
 approach is focused on the construction of an algorithm which
 faces the problem symbolically. What this algorithm is essentially doing is that
 helps us to rewrite $p$ as follows:
 \[ p=c_1
(x,y',\ldots,y^{(k-1)},W_{i,\sigma,\varphi})
[W_{0,k}+W_{1,k}y+W_{2,k}y'+\cdots+y^{(k)}]^{j_{0,1}}\]
\[[[W_{0,k+1,1}+W_{1,k+1,1}y+W_{2,k+1,1}y'+\cdots+y^{(k)}]']^{j_{1,1}}\]
\[\cdots [[W_{0,n,1}+W_{1,n,1}y+W_{2,n,1}y'+\cdots+y^{(k)}]^{(n-k)}]^{j_{n-k,1}}\]
\[+\cdots+c_\nu
(x,y',\ldots,y^{(k-1)},W_{i,\sigma,\varphi})\]
\[[W_{0,k}+W_{1,k}y+W_{2,k}y'+\cdots+y^{(k)}]^{j_{0,1}}\]
\[[[W_{0,k+1,\nu}+W_{1,k+1,\nu}y+W_{2,k+1,\nu}y'+\cdots+y^{(k)}]']^{j_{1,\nu}}\]
\[\cdots [[W_{0,n,\nu}+W_{1,n,\nu}y+W_{2,n,\nu}y'+\cdots+y^{(k)}]^{(n-k)}]^{j_{n-k,\nu}}
+R\]
Where $k$ is given, $n$ is the order of $p$, and $j_{a,b}$ are specific
positive whole numbers. The quantities $W_{i,\sigma,\varphi}$
are undetermined parameters which can take certain values, 
$c_j$ are the coefficients, depending
from the parameters $W_{i,\sigma,\varphi}$ and $R$
a polynomial
of the variables $x,y,y', \ldots, y^{(k-1)}$ and the parameters
$W_{i,\sigma,\varphi}$ called the remainder. Afterwards, we seek for those values
of the parameters which eliminate the remainder. If this is possible, 
then the linear differential polynomial
$W_{0,k}+W_{1,k}y+ \cdots + y^{(k-1)}$, with $W_{0,k},W_{1,k}, \ldots$
evaluated over those values which vanish the remainder, is a factor of
$p$, (where the operation of differentiation has been taken under
consideration). This means that any solution
of $W_{0,k}+W_{1,k}y+ \cdots + y^{(k-1)}=0$ is a solution of the 
equation (\ref{arxh}), too. Since we do it for every $k=0, \ldots,n$,
we collect likewise, all the " linear " solutions of the equation
$p=0$. 
\par
Our method is an extension of a similar procedure, introduced 
by the author and applied in the study of difference equations
and feedback design \cite{kn:difeq},\cite{kn:ieee}.
Its main merit is its computational orientation.
It turns to be a useful tool, implemented on a computer machine and 
gives useful results. 
Moreover, despite our method resembles with the 
approach
of differential algebra, \cite{kn:ritt}, there are some differences. Concretely,
(i) We are working with a single differential polynomial whilst Ritt's 
algorithm deals with sets of differential polynomials. 
(ii) The existence
of the parameters $W_{i,\sigma,\varphi}$ permits us to 
find classes of linear solutions. We can then 
select among them, these particular solutions which satisfy
additional conditions.
For instance, we can search for those values of $W_{i,\sigma,\varphi}$, 
if any, which do not only eliminate the remainder but also 
yield stable linear factors. 
(iii) In the classical Ritt's approach we find the minimum number of 
differential polynomials which generate a differential ideal. In our 
method we check if a given polynomial belongs to a differential ideal, 
produced by \underline{linear} differential polynomials. 
 Throughout the text, ${\bf R}$ and ${\bf Z}^+$ will denote the sets of real
numbers and positive integers, correspondingly.

\section{The Algebraic Framework}
 Let  
${\bf R}[x]$ be the ring of polynomials of a single variable
 with real coefficients. This
polynomial ring is a differential ring too, with the usual 
derivation \cite{kn:ritt},\cite{kn:kapla}. Let $y(x)$ be 
a real function and
$y^{(i)}(x)$, $i=0,1,2,\ldots$ its derivatives. A differential 
polynomial $p$, in $y(x)$ or shortly in $y$, is a polynomial in $y$ and its derivatives
with coefficients in ${\bf R}[x]$. $p$ can be written as follows:
\[ p=\sum_{\lambda=1}^{\varphi}s_{\lambda}x^{a_\lambda}
\prod_{i=0}^n[y^{(i)}(x)]^{\theta_{i,\lambda}}\]
where $s_\lambda \in {\bf R}$ and some of the 
exponents $a_\lambda,\theta_{i,\lambda} \in {\bf Z}^{+} $ are
not equal to zero. The number $n$, which represents the highest 
order derivative of $y(x)$, is called the order of $p$.
An equation of the form $p=0$, with $y(x)$ as unknown function,
is called a polynomial differential equation. Any function which satisfies
it, is called a solution or a general solution. An expression of the form:
$L=\sum_{i=0}^na_iy^{(i)}(x)$, $a_i \in {\bf R}$, is called a
linear differential polynomial and the equation $L=0$ 
a linear 
differential equation. Its solutions are called 
" linear " solutions of order $n$. 
\par
\noindent
The highest derivative of $y(x)$, appeared in the polynomial
$p$, is called the {\it leader}. Let 
$p_1=s_1x^{a_1}\prod_{i=0}^n[y^{(i)}(x)]^{\theta_{i,1}}$,
$p_2=s_2x^{a_2}\prod_{i=0}^n[y^{(i)}(x)]^{\theta_{i,2}}$ 
be two, not identical, terms of $p$. This means that there is
at least one index $k$, $1\le k \le n$, such that $\theta_{k,1} \ne \theta_{k,2}$
or $\theta_{i,1} = \theta_{i,2}$, $i=0,1,\ldots n$ and
$a_1 \ne a_2$. We say that the term $p_2$ is ordered higher 
than $p_1$ with respect to 
{\it lexicographical } order and we write $p_1\prec p_2$, 
if either there is an index
$s$ such that $\theta_{s,1}<\theta_{s,2}$ and
$\theta_{j,1}=\theta_{j,2}$, $j=s+1, \ldots, n$
or
$\theta_{i,1}=\theta_{i,2}$, $i=0,1,\ldots,n$ and $a_1<a_2$.
By means of this rank we can order all the terms of $p$ in
an ascending way. The term which is 
ordered higher, is called the {\it maximimum} term of $p$.
\par
\noindent
The {\it differential ideal}, generated
by a finite set of differential polynomials: \newline 
${\bf \Phi}=
\{\Phi_1,\Phi_2,\ldots,\Phi_m\}$ and denoted by
$[{\bf \Phi}]$ is a set which consists of all differential 
polynomials that can be formed of elements in ${\bf \Phi}$
 by multiplication with arbitrary polynomials, addition and differentiation.
\par
\noindent
Let ${\cal W}=\{W_{i,\sigma,\varphi}\}$ be a set of undetermined 
parameters, taking values in ${\bf R}$. A {\it Formal-k-Factorization} of
$p$, 
denoted by $Formal(p,k)$, is an expression of $p$ of the form:

\[ Formal(p,k)=\sum_{\mu=1}^{\nu}c_\mu
(x,y'(x),\ldots,y^{(k-1)}(x),W_{i,\sigma,\varphi})
\cdot \]
\begin{equation}\label{ff}
\cdot \prod_{i=0}^{n-k}[(W_{0,k+i,\mu}+W_{1,k+i,\mu}y(x)+
W_{2,k+i,\mu}y'(x)+\cdots+y^{(k)}(x))^{(i)}]^{j_i}+R
\end{equation}
where the coefficients 
$c_\mu
(x,y'(x),\ldots,y^{(k-1)}(x),W_{i,\sigma,\varphi})$ and the remainder
$R$ are polynomials of
the terms $x,y'(x),\ldots,y^{(k-1)}(x),W_{i,\sigma,\varphi}$
only. Some of the exponents $j_i \in {\bf Z}^+$ may be equal to zero.
Sometimes, (\ref{ff}) is written briefly as

\[Formal(p,k)=\sum_{\mu=1}^{\nu}c_{\mu}\prod_{i=0}^{n-k}[L_{i,\mu}^{(i)}]^{j_{i}}+R\]
\noindent
where we used the notation $L_{c,k}$ for the linear differential
polynomial $W_{0,k+i,\mu}+W_{1,k+i,\mu}y(x)+
W_{2,k+i,\mu}y'(x)+\cdots+y^{(k)}(x)$. 
\par
\noindent
We can take different expressions of the $Formal(p,k)$ of a 
concrete differential polynomial $p$, 
by giving  to the parameters $W_{i,\sigma,\varphi}$ 
certain values. Such procedures are called
{\it evaluations} of the $Formal(p,k)$. A most rigorous approach is the following:
Let ${\cal W}=\{W_{i,\sigma, \varphi}\}$ be the set of the variables, appeared in the
Formal - k - Factorization of a given polynomial $p$.
By arranging the parameters in an increasing order we form the vector
${\cal W}=(W_{i_h,\sigma_h,\varphi_h})_{h=1,2, \ldots,n}$.
Let ${\bf r}=(a_h)_{h=1,2,\ldots,n}$ be a vector of real numbers,
 which has the same length
with the vector ${\cal W}$. We say that the parameters ${\cal W}$ follow the rules
${\bf r}$ and we write ${\cal W} \to  {\bf r}$ 
if the following substitution is valid:
$W_{i_h,\sigma_h,\varphi_h}=a_h$, $h=1,2,\ldots,n$.
Let $M$ a set of rules, $M=\{{\bf r}_1,{\bf r}_2, \ldots, {\bf r}_\lambda\}$ then
\[ \left. \begin{array}
{c}
Formal(p,k)\\
\end{array} \right|_{M} =\bigcup_{\nu=1}^\lambda 
\{ \sum_{\mu=1}^{\nu}c_\mu\cdot \prod_{i=0}^{n-k}[(W_{0,k+i,\mu}+W_{1,k+i,\mu}y(x)+\]
\[+W_{2,k+i,\mu}y'(x)+
\cdots+y^{(k)}(x))^{(i)}]^{j_i}+R
 \quad \hbox{with} 
\quad {\cal W} \to {\bf r}_\nu \in M \}\]
The set of substitutions $M$, may be finite or infinite. 
When we evaluate $Formal(p,k)$ over $M$, the linear 
differential polynomials $L_{i,\mu}$ and the remainder $R$
take specific values, these are denoted
by
$ \left. \begin{array}
{c}
L_{i,\mu}\\
\end{array} \right|_{M}$ and
$\left. \begin{array}
{c}
R\\
\end{array} \right|_{M}$.
A case of particular interest is when we can find values
of the parameters which eliminate the remainder $R$. 
Whenever this happens, $p$ is a " combination " of linear
differential polynomials or, in a more formal language, $p$ 
is a member of the differential 
ideal produced by these linear differential polynomials.
Relevant is the following theorem:

\smallskip
THEOREM 1. Let $p$ be a differential polynomial of order $n$,
let $k$ be given and $Formal(p,k)=$
$\sum_{\mu=1}^{\nu}c_{\mu}
\prod_{i=0}^{n-k}[L_{i,\mu}^{(i)}]^{j_{i}}+R$ its Formal-k-Factorization.
Let us suppose that there is a set of rules, denoted by ${\cal R}$, 
which eliminates the remainder $R$, i.e. 
$ \left. \begin{array}
{c}
R\\
\end{array} \right|_{\cal R}=0$, 
then $p \in [\left. \begin{array}
{c}
L_{i,\mu}\\
\end{array} \right|_{\cal R},$$i=0,\ldots,n-k,\mu=1,\ldots,\nu]$.
\smallskip

PROOF. The proof comes straightforward from the definition
of the Formal-k-Factorization.
$ \square$
\par
\noindent
It is obvious that if the linear differential equations
$\left. \begin{array}
{c}
L_{i,\mu}\\
\end{array} \right|_{\cal R}=0$
have a common solution, this is a solution of the nonlinear
equation $p=0$, too. This is the cornerstone of our approach.
\smallskip

EXAMPLE 1.
Let us consider the differential polynomial
\[p=4y'''-4(y'')^2+y'y''-\frac{1}{16}(y')^2-1\]
 then,
\[Formal(p,2)=4(W_{0,3,1}+W_{1,3,1}y+W_{2,3,1}y'+y'')'-4(W_{0,2,2}+W_{1,2,2}y+W_{2,2,2}y'+y'')^2\]
\[+[8W_{0,2,2}-4W_{2,3,1}+8W_{1,2,2}+(8W_{2,2,2}+1)y']\cdot (W_{0,2,3}+W_{1,2,3}y+W_{2,2,3}y'+y'')+R\]
where the remainder
$R$ is
$R=4W_{0,2,2}^2$$-8W_{0,2,3}W_{0,2,2}$
$\newline +4W_{0,2,3}W_{2,3,1}$$-1$
$+(8W_{1,2,2}W_{0,2,2}$$-8W_{1,2,3}W_{0,2,2}$$-8W_{0,2,3}W_{1,2,2}$$+4W_{1,2,3}W_{2,3,1})y$
$\newline +(4W_{1,2,2}^2-8W_{1,2,2}W_{1,2,3})y^2$
$+(8W_{2,2,2}W_{0,2,2}$$-8W_{2,2,3}W_{0,2,2}$
$-W_{0,2,3}-4W_{1,3,1}$$-8W_{0,2,3}W_{2,2,2}$$\newline 
+4W_{2,2,3}W_{2,3,1})y'$
$ +(4W_{2,2,2}^2$$-8W_{2,2,2}W_{2,2,3}$$-W_{2,2,3}-\frac{1}{16})(y')^2$
$+(-W_{1,2,3}+8W_{1,2,2}W_{2,2,2}$
$-8W_{1,2,3}W_{2,2,2}$
$-8W_{1,2,2}W_{2,2,3})yy'$. The following rules eliminate
the remainder:
${\bf r}_1=\{ W_{0,3,1}=s,W_{0,2,2}=\omega, $
$W_{2,2,3}=\varphi,$
$W_{0,2,3}=k,$
$W_{1,3,1}=\frac{1}{4}\omega-\frac{(4\omega^2-1)\varphi}{4k},$
$W_{1,2,3}=0,$
$W_{2,2,2}=-\frac{1}{8},$
$W_{1,2,2}=0,$
$W_{2,3,1}=2\omega+\frac{1-4\omega^2}{4k}\}$
and
${\bf r}_2=\{ W_{0,3,1}=s,W_{0,2,2}=\omega, $
$W_{2,2,3}=\varphi,$
$W_{0,2,3}=k,$
$W_{1,3,1}=\frac{1}{4k}[\varphi+(\omega-2k)(k-4(\omega-2k)\varphi)],$
$W_{1,2,3}=0,$
$W_{2,2,2}=2\varphi+\frac{1}{8},$
$W_{1,2,2}=0,$
$W_{2,3,1}=2\omega+\frac{1-4\omega^2}{4k}\}$
with $\omega,\varphi,k,s \in {\bf R}$.
Indeed, for instance
\[\left. \begin{array}
{c}
Formal(p,2)\\
\end{array} \right|_{{\bf r}_1}=4(s-\frac{1}{4}\left(\omega+\frac{(4\omega^2-1)\varphi}{k}\right)y+\left(2\omega+\frac{1-4\omega^2}{4k}\right)y'+y'')'\]
\[-4(\omega-\frac{y'}{8}+y'')^2+\left(\frac{4\omega^2}{k}-\frac{1}{k}\right)
\cdot (\omega +\varphi y'+y'')\]
and the differential ideal which
contains is 
$[s-\frac{1}{4}\left(\omega+\frac{(4\omega^2-1)\varphi}{k}\right)y,$
$+\left(2\omega+\frac{1-4\omega^2}{4k}\right)y'+y'',$
$\omega-\frac{y'}{8}+y'',$
$\omega +\varphi y'+y'']$
$\omega,\varphi,s,k,\in {\bf R}$.
By setting $\omega=\varphi=s=0,k=1$ we take
the simplified ideal
$[y'',-\frac{y'}{8}+y'']$.
We can take analogous expression for
$\left. \begin{array}
{c}
Formal(p,2)\\
\end{array} \right|_{{\bf r}_2}$.

\section{Detection of the Linear Solutions}
The scope of this section is to present the algorithm which 
constructs for a given $k$, a special $Formal(p,k)$ with a linear differential polynomial
 as a common factor to each term. We denote this factor by $L_{c,k}$. Afterwards, by finding
proper sets of values for the parameters $W_{i,\sigma,\varphi}$, 
we eliminate the remainder. It is then clear that any solution
of the linear equation $L_{c,k}=0$, where the polynomial $L_{c,k}$
has been evaluated over this set, is a solution of the original system, too.
By repeating the whole procedure for every $k=0,\ldots,n$ 
we discover all the linear solutions. As we pointed out, the crucial
issue is how can we eliminate the remainder. This is carried out 
by solving a system of algebraic equations. 
Finally, we have to elucidate that in this paper we 
 do not take into account initial conditions. We are only focused on how we obtain
 general solutions, that is solutions which " contain " constants.
We present now 
the algorithm upon discussion.
\par
\noindent
Let us suppose that an algorithm which solves an algebraic
system of polynomial equations, is available. These 
algorithms are classical in computational algebra and there
are many of them in the literature \cite{kn:cox1}. We 
name a such algorithm as $SysAlgEqs$.
\par
\vskip 10 pt 
\noindent 
\underline{\bf THE DIF-FORMAL ALGORITHM} 
\vskip 10 pt 
\small {\sf
\par
\noindent {\bf Input:} \begin{itemize}
\item A differential polynomial $p$ of order $n$.
\item A set of undetermined parameters
${\cal W}=\{ W_{i,\sigma,\varphi}\}$, taking values in ${\bf R}$.
\end{itemize}
\vskip 5 pt
\noindent
{\bf Output:}
The quantities $S_k$, $k=0,\ldots,n$
}
\vskip 5 pt
\noindent {\bf FOR} $k=0$ {\bf TO} $n$ 
\par
\begin{description}
\item{\bf Step 1:} $R=p$, $\mu=0$.

\item{\bf Step 2:} {\bf REPEAT} the 
following steps {\bf UNTIL} $R$ does not contain terms
of order $\ge k$.
\begin{description}
\item{\bf Step 2a:} Set $\mu=\mu+1$,
\item{\bf Step 2b:} Find the maximum term of $R$, 
with respect to the lexicographical order.
We denote it by
\[ p_\mu=s_\mu \cdot x^{a_\mu}\cdot 
\prod_{i=0}^n[y^{(i)}(x)]^{\lambda_{i,\mu}} \]
where $a_\mu,\lambda_{i,\mu}$ are positive integers and
$y^{(i)}(x)$ the derivatives of $y(x)$ of order $i$.
At the first iteration $s_\mu$ is a constant, then it becomes
a function of the free parameters  $W_{i,\sigma,\varphi}$ as well.
\item{\bf Step 2c:} Construct the linear formal differential polynomials:
\[ L_{c,k}=W_{0,k}+W_{1,k}y(x)+W_{2,k}y'(x)+ \cdots +W_{k,k}y^{(k-1)}(x)+y^{(k)}(x)= \]
\[ =W_{0,k}+\sum_{j=0}^{k-1}W_{j+1,k}y^{(j)}(x)+y^{(k)}(x) \]
\[ L_{i,k}=W_{0,k+i,\mu}+\sum_{j=0}^{k-1}W_{j+1,k+i,\mu}y^{(j)}(x)+y^{(k)}(x), \quad i=1, \ldots,n-k \]
\item{\bf Step 2d:} Execute the operation:
\[ R=R - s_\mu \cdot x^{a_\mu}\cdot \prod_{i=0}^{k-1} [y^{(i)}(x)]^{\lambda_{i,\mu}} 
\cdot [L_{c,k}^{(0)}]^{\lambda_{k,\mu}}
\cdot \prod_{i=1}^{n-k} [L_{i,k}^{(i)}]^{\lambda_{i,\mu}} \]
\item{\bf END of REPEAT}
\end{description}
\item{\bf Step 3:} By means of the $SysAlgEqs$-Algorithm
we find the set $S_k$ 
of those values of the parameters, $W_{i,\sigma,\varphi}$, 
which eliminate the remainder. In other words:
$\left. \begin{array}
{c}
R\\
\end{array} \right|_{S_k}=0$.

\item{\bf END of FOR}
 \end{description}
 \normalsize 
\noindent

It is obvious that the DIF-FORMAL algorithm terminates 
after a finite number of iterations. 

THEOREM 2. Let $S_k$, $k=0,\ldots,n$ be the outputs of the DIF-FORMAL
Algorithm. If $S_k \ne \emptyset$ for some values of $k$, then the 
solutions of the linear differential equations
$\left. \begin{array}
{c}
L_{c,k}\\
\end{array} \right|_{S_k}=0$
are solutions of the nonlinear polynomial differential equation
$p=0$, too. 
\smallskip

PROOF. Let $p$ a differential polynomial and $k$ fixed.
By substituting backwards the successive results of the 
step 2d
we find that
\[ p=\sum_{\mu=1}^m s_{\mu}x^{a_\mu}\prod_{i=0}^{k-1}[y^{(i)}]^{\lambda_{i,\mu}}
\cdot [L_{c,k}^{(0)}]^{\lambda_{i,\mu}}\cdot \prod_{i=1}^{n-k}[L_{i,k}^{(i)}]^{\lambda_{i,\mu}}\]
This is the Formal-k-Factorization of $p$
with $c_\mu=s_{\mu}x^{a_\mu}\prod_{i=0}^{k-1}[y^{(i)}]^{\lambda_{i,\mu}}$
and $L_{c,k}$ as a common factor in all the terms
but the remainder. We denote it by $CFormal(p,k)$.
Let us now suppose that $S_k \ne \emptyset$, this means that 
$\left. \begin{array}
{c}
R\\
\end{array} \right|_{S_k}=0$ and thus, the linear differential polynomial
$L_{c,k}$, evaluated over $S_k$, is a factor of every term
of $\left. \begin{array}
{c}
CFormal(p,k)\\
\end{array} \right|_{S_k}=0$. This implies that 
any solution of the linear differential equation
$\left. \begin{array}
{c}
L_{c,k}\\
\end{array} \right|_{S_k}=0$
is also a solution of the nonlinear equation $p=0$. Since this
argument is true for any $k$, the theorem has been 
proved.$\square$
\par
\noindent
The above result can be restated, using ideals, in the following
way.

\smallskip

COROLLARY 1. Let $S_k$, $k=0,\ldots,n$ be the outputs of the DIF-FORMAL
Algorithm. If $S_k \ne \emptyset$ for some values of $k$, then
$p \in [ \left. \begin{array}
{c}
L_{c,k}\\
\end{array} \right|_{S_k}]$, $k=0,\ldots,n$. 
\smallskip

In order to eliminate the remainder we have to solve a system of algebraic polynomial equations. This 
can be done via several methods. Groebner basis, \cite{kn:cox1}, is a popular powerfull tool, with
satisfactory results.
\par
We can extend the algorithm posted above
toward different directions.
(a) Instead of the constant term in the common linear 
differential polynomial we can use a single polynomial of $x$ 
of a given degree with parametrical coefficients. These 
coefficients can be determined through the remainder elimination, too.
(b) Instead of linear differential polynomial, with constant
coefficients, we can use linear differential polynomials 
with variable coefficients which correspond to
solvable differential equations (Euler equations, for instance).
The last two cases are currently under studying.
\par \noindent
\smallskip
EXAMPLE 2. Let us consider the differential equation
\[(y'')^3-2y'(y'')^2-4(y'')^2+(y')^2y'''+4y'y'''+4y'''=0\]
 or $p=0$, where $y=y(x)$ is the unknown function
 and the order of $p$ is $n=3$. We want to find all the "linear"
 solutions which are included into this equation. The application
 of the DIF-FORMAL algorithm gave the following results:
 $y(x)=c, c \in {\bf R}$, $y(x)=c_1+xc_2$, $y(x)=c_1e^x-2x+c_2$,
 $c_1,c_2 \in {\bf R}$. To clarify our ideas and to indicate how
 the algorithm works in practice, we shall present the case 
 $k=2$ in details.
 Since $k=2$ we are going to detect linear polynomials of second
 order included into the original equation. At the first iteration
 the maximum term is $[(y')^2]y'''$. We construct 
 the linear polynomial:
 \[ L_{c,2}=W_{0,2}+W_{1,2}y+W_{2,2}y'+y'' \]
 and we execute the subtraction:
 $p_1=p-[(y')^2]\cdot L_{c,2}'$. This operation will eliminate
 the $y'''$ term. In the next iteration $4y'y'''$ is the maximum
 term. Since we are looking for a common second order linear equation
 included into $p$, we use the same $L_{c,2}$ and 
 we calculate $p_2=p_1-4y'\cdot L_{c,2}$. Working this way
 we finally get:
 \[ Formal(p,2)=[(y')^2+4y'+4]\cdot L_{c,2}'+L_{c,2}^3+
 [-3W_{0,2}-3W_{1,2}y-\]
 \[-3W_{2,2}y'-2y'-4] L_{c,2}^2+[3W_{0,2}^2+6W_{1,2}W_{0,2}y+
 6W_{2,2}W_{0,2}y'+4W_{0,2}y'+\]
 \[+8W_{0,2}+3W_{1,2}^2y^2+3W_{2,2}^2(y')^2+3W_{2,2}(y')^2
 -4W_{2,2}+8W_{1,2}y+4W_{2,2}y'+\]
 \[+4W_{1,2}yy'+6W_{1,2}W_{2,2}yy']L_{c,2}+R\]
 where $R=(-3W_{1,2}W_{0,2}^2$
 $-8W_{1,2}W_{0,2}$
 $+4W_{1,2}W_{2,2})y$
 $+(-3W_{0,2}W_{1,2}^2-4W_{1,2}^2)y^2$
 $-W_{1,2}^3y^3$
 $+(-3W_{2,2}W_{0,2}^2$
 $-2W_{0,2}^2$
 $-4W_{2,2}W_{0,2}$
 $+4W_{2,2}^2$
 $-4W_{1,2})y'$
 $+(-4W_{0,2}W_{1,2}$
 $\newline -6W_{0,2}W_{2,2}W_{1,2}$
 $-4W_{2,2}W_{1,2})yy'$
 $+(-3W_{2,2}W_{1,2}^2$
 $-2W_{1,2}^2)y^2y'$
 $+(-3W_{0,2}W_{2,2}^2$
 $-3W_{0,2}W_{2,2}$
 $-4W_{1,2})(y')^2$
 $+(-3W_{1,2}W_{2,2}^2$
 $-3W_{1,2}W_{2,2})y(y')^2$
 $+(-W_{2,2}^3-W_{2,2}^2$
 $-W_{1,2})(y')^3$
 $+4W_{2,2}W_{0,2}$
 $-W_{0,2}^3-4W_{0,2}^2$.
 The values $(W_{0,2}=-2,W_{1,2}=0,W_{2,2}=-1)$ and
 $(W_{0,2}=0,W_{1,2}=0,W_{2,2}=0)$, eliminate the remainder.
 By substituting them to $L_{c,2}$ we get the linear differential
 equations $y''-y'-2=0$, $y''=0$. The general solutions of those
 equations will be solutions of the original nonlinear differential
 equation too. Therefore $y(x)=c_1e^x+c_2-2x$, $y(x)=c_2x+c_1$ are solutions 
 of the equation $p=0$, too. Working similarly, we obtain for 
 $k=0$ and $k=1$ analogous results. What we have actually proved
 is that $p$ is a member of the differential ideal $[y''-y-2,y'']$.
 \par
 \noindent
\smallskip
EXAMPLE 3. Let us illustrate now the case where, instead 
 of constant terms, we have polynomials of $x$. We consider the
 differential equation:
 \[(x^2-x)y'+xy''-x^2y'''+(-2x^3-3x^2+3x)=0\]
 In this case we work with the extended
 algorithm which uses common linear factors of the form:
 $y'+W_{1,1}y+W_{-2,0,1}x^2+W_{-1,0,1}x+W_{0,0,1}$.
 This modified DIF-FORMAL Algorithm will give that 
 $W_{1,1}=-1$, $W_{-2,0,1}=1$, $W_{-1,0,1}=3$ and $W_{0,0,1}$
 free or
 $W_{1,1}=0$, $W_{-2,0,1}=0$, $W_{-1,0,1}=-2$ and 
 $W_{0,0,1}=-5$. The corresponding general linear solution is
 $y(x)=c_1e^x+x^2+5x+c_2$.


\begin{thebibliography}{99}

\bibitem{kn:difeq} S. Kotsios,
" On Detecting Solutions of Polynomials Nonlinear Difference 
Equations ".
{\it Journal of Difference Equations and Applications, Vol. 8 (6), pp 551-571, 2002}.

\bibitem{kn:kapla} I. Kaplansky, "Differential Algebra", Herman
Paris, 1957.
\bibitem{kn:cox1} D. Cox, J. Little, D. O'Shea. (1997). "Ideals, Varieties
and Algorithms". {\it Springer-Verlag, New York}.
\bibitem{kn:cheng1} D. Behloul, S.S. Cheng (2005) " Computation
of All Polynomial Solutions of a Class of Nonlinear Differential Equations "
to appear.
 \bibitem{kn:ritt} J.F.Ritt (1950), "Differential Algebra", {\it American Mathematical
Society, Providence, RI.}
 \bibitem{kn:cheng2} W. R. Li, S.S. Cheng, T.T. Lu, 
 " Closed form solutions of iterative 
 functional differential equations ", 
 {\it Appl. Math. E-Notes, 1 (2001),1-4}
 
\bibitem{kn:glad1} S.T. Glad, " Differential Algebraic Modeling of Non-Linear
Systems ", {\it In proceedings of MTNS-89, Amsterdam, 97-105, 1989.}


\bibitem{kn:cheb} E.S. Cheb-Terrab, L.G.S. Duarte,
 L.A.C.P. da Mota, 
 " Computer algebra solving of first order ODEs 
 using symmetry methods. " 
 {\it Comput. Phys. Comm. 101 (1997), No 3, 254-268.}
 

\bibitem{kn:ieee} S. Kotsios, " A Symbolic Computational 
Algorithm for Designing Feedback Stabilizers
of Nonlinear Systems by Checking Positivity of
 Polynomials ", It has been submitted to IEEE TAC as technical note. 


\end{thebibliography}
\end{document}